# On a Weakened Form of Polignac Conjecture

Shaohua Zhang [1, 2]

[1] School of mathematics of Shandong University, People's Republic of China 250100

[2] The key lab of cryptography technology and information security, Ministry of Education,

Shandong University, People's Republic of China 250100

E-mail: shaohuazhang@mail.sdu.edu.cn

**Abstract:** Polignac [1] conjectured that for every even natural number $2k(k \geq 1)$, there exist infinitely many consecutive primes $p_n, p_{n+1}$ such that $p_{n+1} - p_n = 2k$. A weakened form of this conjecture states that for every $k \geq 1$, there exist infinitely many primes $p$, $q$ such that $p - q = 2k$. Clearly, the weakened form of Polignac's conjecture implies that there exists an infinite sequence of positive integers $x_1, ..., x_m, ...$ such that $x_1(2k + x_1), ..., x_m(2k + x_m), ...$ are pairwise relatively prime. In this note, we obtain a slightly stronger result than this necessary condition. This enables us to find a common property on some special kinds of number-theoretic functions (such as $2^x - 1$) which likely represent infinitely many primes by rich literatures and a lot of research reports.

However, the function $2^{2^x} + 1$ does not have this property. Does it imply that the number of Fermat primes is finite? Hardy and Wright [7] conjectured that the number of Fermat primes is finite. Nevertheless, they did not give any reasons and explanations. By factoring Fermat number, many people believe that the conjecture in [7] holds. Does our work explain this phenomenon? We will consider further this problem in another paper.

Based on our work, one could give a new sufficient condition that there are an infinite number of twin primes (Sophie-Germain primes or Mersenne primes).

**1 INTRODUCTION** In 1849, Polignac [1] conjectured that for every even natural number $2k(k \geq 1)$, there exist infinitely many consecutive primes $p_n, p_{n+1}$ such that $p_{n+1} - p_n = 2k$. Paulo Ribenboim [2] pointed out that the weakened form of this conjecture has never been proven: For every $k \geq 1$, there exist primes $p$ and $q$ such that $p - q = 2k$, let alone there exist infinitely many primes $p$ and $q$ such that $p - q = 2k$, as we know, even if $k = 1$, it is just the famous twin prime conjecture which is still open. Jingrun Chen proved that there are infinitely many numbers $x$ such that $x + 2$ has at most two prime divisors [**19**]. Amazingly enough, this weakened form is analogous to Goldbach's conjecture which also is open and states that every even integer $2k \geq 4$ is the sum of two primes. Jingrun Chen's well known theorem states every sufficiently large even number can be written as the sum of a prime and the product of at most two prime divisors [**19**]. Combining with these two problems together, can you prove that every even integer $2k \geq 4$ is the sum or difference of two primes? This becomes an interesting problem. But, we do not know whether it is proved or not. It seems that these conjectures should hold. Particularly, for every even natural number $2k$, not only there are primes $p$ and $q$ such that $p - q = 2k$, but also such pairs of $p$ and $q$ should be infinitely many. In 1923, Hardy and Littlewood [**3**: conjecture B] conjectured that for every $k \geq 1$, the number of $n > 1$ such that $p_n \leq x$ and $p_{n+1} = 2k + p_n$ is about $(C + o(1)) \times (\prod_{\substack{q|k \\ odd \\ prime}} \frac{q-1}{q-2}) \times \frac{x}{(\log x)^2}$. Thus, Hardy-Littlewood conjecture implies the following weakened form of Polignac's conjecture holds.

**The Weakened Form of Polignac's Conjecture:** For every $k \geq 1$, there exist infinitely many primes $p$, $q$ such that $p - q = 2k$.

Recently, by reading Pintz's paper "Landau's problems on primes", we also know that Kronecker mentioned this problem in 1901.

This weakened form of Polignac's conjecture can be viewed as the special case of Dickson's conjecture [**4**] or Schinzel-Sierpinski's Conjecture [**5**]. By Schinzel-Sierpinski's Conjecture, the necessary condition that several irreducible univariable polynomials $f_1(x), \cdots, f_s(x)$ with integer coefficients represent simultaneously primes for infinitely many $x$ is of that their leading coefficients are positive, and, there does not exist any integer $n > 1$ dividing all the products $\prod_{i=1}^{s} f_i(k)$, for every integer $k$. In this note, one of our aims is to give another necessary condition, although precise conjectures do not seem to have been formulated in the literature for that several multivariable number-theoretic functions represent simultaneously primes for infinitely many integral points. Note that pairwise distinct primes are pairwise relatively prime. Therefore, if $f_1(x), \cdots, f_s(x)$ represent simultaneously primes for infinitely many integers $x$, then there exists an infinite sequence of integers $x_1, \ldots, x_k, \ldots$ such that $\prod_{j=1}^{s} f_j(x_1), \ldots, \prod_{j=1}^{s} f_j(x_k), \ldots$ are pairwise relatively prime. In fact, it is not difficult to prove that these two necessary conditions are equivalent when $f_1(x), \cdots, f_s(x)$ are irreducible polynomials. Thus, Schinzel-Sierpinski's Conjecture can be re-stated as follows: if $f_1(x), \cdots, f_s(x)$ are irreducible polynomials, and there exists an infinite sequence of integers $x_1, \ldots, x_k, \ldots$ such that $\prod_{j=1}^{s} f_j(x_1), \ldots, \prod_{j=1}^{s} f_j(x_k), \ldots$ are pairwise relatively prime, then $f_1(x), \cdots, f_s(x)$ represent simultaneously primes for infinitely many integers $x$.

More generally, let $f_1(x_1, \cdots, x_k), \ldots, f_s(x_1, \cdots, x_k)$ be $s$ multivariable number-theoretic functions. In this note, we restrict that a multivariable number-theoretic function is a map from $N^k$ to $Z$. If $f_1(x_1, \cdots, x_k), \ldots,$ and $f_s(x_1, \cdots, x_k)$ represent simultaneously primes for infinitely many integral points $x = (x_1, \cdots, x_k)$, then there exists an infinite sequence of integral points $(x_{11}, \cdots, x_{k1}), \ldots, (x_{1i}, \cdots, x_{ki}), \ldots$ such that $\prod_{j=1}^{s} f_j(x_{11}, \ldots, x_{k1}), \ldots, \prod_{j=1}^{s} f_j(x_{1i}, \ldots, x_{ki}), \ldots$ are pairwise relatively prime. This should be viewed as a natural necessary condition that several multivariable number-theoretic functions represent simultaneously primes for infinitely many integral points. Particularly, the weakened form of Polignac's conjecture is equivalent to the following:

The number-theoretic functions $f_1(x) = x$ and $f_2(x) = x + 2k$ represent simultaneously primes for infinitely many variable $x$.

Thus, if the weakened form of Polignac's conjecture holds, then there is an infinite sequence of positive integers $x_1, \ldots, x_m, \ldots$ such that $x_1(2k + x_1), \ldots, x_m(2k + x_m), \ldots$ are pairwise relatively prime. In this note, we obtain a slightly stronger result than this necessary condition as follows:

**Theorem 1:** For every integer $a \geq 1$, there is a constant $c$, such that provided $n > c$, there is a positive integer $x$ such that $x > 1$, $x \in Z_n^*$ and $x + 2a \in Z_n^*$, where $Z_n^* = \{x \in N \mid 1 \leq x < n, (x, n) = 1\}$.

Namely, let $a$ be a natural number. Then there is a constant $c$ (only depending on $a$), so that for integers $n > c$, there is an integer $x$ with $1 < x < x + 2a < n$, moreover, $x$ and $x + 2a$ relatively prime to $n$. Clearly, this implies that there exists an infinite sequence of positive integers $x_1, \ldots, x_m, \ldots$ such that $x_1(2k+x_1),...,x_m(2k+x_m),...$ are pairwise relatively prime. Based on this result, we will try to give a new sufficient condition that there are an infinite number of twin primes. We have:

**Conjecture 1:** Let $n>3$ be a positive integer. If $x>1$ is the smallest integer such that $x \in Z^*_{n!}$ and $x+2 \in Z^*_{n!}$, then $x$ and $x+2$ are primes.

## 2 THE PROOF OF THEOREM 1
We now prove the theorem. Firstly, we need the following lemmas.

**Lemma 1:** Let $p$ be an odd prime number, and let $(a,b)=1$. If $(ab,p)=p$, then we have

$$\#\{x \mid x \in Z^*_p, (ax+b,p)=1\} = p-1 ; \text{ Otherwise, } \#\{x \mid x \in Z^*_p, (ax+b,p)=1\} = p-2.$$

**Proof:** Easy.

**Corollary 1:** Let $e$ be a positive integer, and let $(a,b)=1$. If $(ab,p)=p$, then

$$\#\{x \mid x \in Z^*_{p^e}, (ax+b, p^e)=1\} = p^{e-1}(p-1) ; \text{ Otherwise}$$

$$\#\{x \mid x \in Z^*_{p^e}, (ax+b, p^e)=1\} = p^{e-1}(p-2)$$

**Lemma 2:** Let $e$ be a positive integer, and let $(a,b)=1$. If $(ab,2)=2$, then we have

$$\#\{x \mid x \in Z^*_{2^e}, (ax+b, 2^e)=1\} = 2^{e-1} ; \text{ Otherwise, } \#\{x \mid x \in Z^*_{2^e}, (ax+b, 2^e)=1\} = 0.$$

**Proof:** Easy.

**Lemma 3:** Denote the standard factorization of $n$ by $n = 2^r p_1^{e_1} \cdots p_k^{e_k}$, where $p_i$ is odd prime. Let $(a,b)=1$ and $2 \mid ab$. If $(ab,n) = 2^l p_1^{f_1} \cdots p_m^{f_m}$ with $k \geq m$ and $r \geq l \geq 1$, then when $k > m$,

$$\#\{x \mid x \in Z^*_n, (ax+b,n)=1\} = 2^{r-1} \prod_{i=1}^{m} p_i^{e_i-1}(p_i-1) \prod_{i=m+1}^{k} p_i^{e_i-1}(p_i-2) = \varphi(n) \prod_{i=m+1}^{k} \frac{p_i-2}{p_i-1}, \text{ and}$$

when $k=m$, $\#\{x \mid x \in Z^*_n, (ax+b,n)=1\} = \varphi(n)$, where $f_i \geq 1$ for $1 \leq i \leq m$.

**Proof:** Set $F_n = \#\{x \mid x \in Z_n^*, (ax+b,n) = 1\}$. Using the Chinese Remainder Theorem, it is easy to prove that $F_n$ is a multiplicative function. Hence Lemma 3 is true by Corollary 1 and Lemma 2.

**Corollary 2:** Let $n = 2^r p_1^{e_1} \cdots p_k^{e_k}$. Let $(a,b) = 1$ and $2 \mid ab$. If $(ab,n) = p_1^{f_1} \cdots p_m^{f_m}$, then $\#\{x \mid x \in Z_n^*, (ax+b,n) = 1\} = p_1^{e_1-1}(p_1-1) \cdots p_m^{e_m-1}(p_m-1) \prod_{i=m+1}^{k} p_i^{e_i-1}(p_i-2)$. Particularly, if $(ab,n) = 1$ and $2 \mid ab$, then, $\#\{x \mid x \in Z_n^*, (ax+b,n) = 1\} = \varphi(n) \prod_{i=1}^{k} \dfrac{p_i - 2}{p_i - 1}$.

**Lemma 4:** Let $a$ be a positive integer. There is a constant $c$ such that $\varphi(n)/2^k > a$ for every $n > c$, where $k$ is the number of the distinct prime factors of $n$.

**Proof:** Let $2^m > a \geq 2^{m-1}$ and $c = (p_1 \cdots p_{\pi(4a+1)})^{m+2}$, where $p_1, \cdots, p_{\pi(4a+1)}$ are all prime numbers less than or equal to $4a+1$. If $n$ has a prime factor $q > 4a+1$, we write $n = q^e(n/q^e)$ with $(q^e, n/q^e) = 1$. So $\varphi(n)/2^k = q^{e-1}(q-1)\varphi(n/q^e)/2^k \geq q^{e-1}(q-1)/4 > a$. Now, suppose that $n$ has no a prime factor $> 4a+1$. Let $n = p_1^{e_1} \cdots p_{\pi(4a+1)}^{e_{\pi(4a+1)}}$, where $e_i \geq 0$. Note that $n > c = (p_1 \cdots p_{\pi(4a+1)})^{m+2}$. Hence there is always $e_i > m+2$ for some $i$. So we have

$$\varphi(n)/2^k = p_i^{e_i-1}(p_i-1)\varphi(n/p_i^{e_i})/2^k \geq p_i^{e_i-1}(p_i-1)/4 > a.$$

**Corollary 3:** Let $a$ be a positive integer. There is a constant $c$ such that $\varphi(n) > a$ for every $n > c$.

**The proof of Theorem 1:** Denote the standard factorization of $n$ by $n = 2^r p_1^{e_1} \cdots p_k^{e_k}$. By Lemma 3, we have $\#\{x \mid x \in Z_n^*, (x+2a,n) = 1\} \geq \varphi(n)/2^k$ when $r \geq 1$. By Corollary 2, we have also $\#\{x \mid x \in Z_n^*, (x+2a,n) = 1\} \geq \varphi(n)/2^k$ when $r = 0$.

Let $\{x \mid x \in Z_n^*, (x+2a,n) = 1\} = \{x_1, \cdots, x_t\}$. Thus, we have $t \geq \varphi(n)/2^k$. Without loss of generality, we assume that $1 \leq x_1 < x_2 < \ldots < x_t$. By Lemma 4, we know that there is a constant $c$ such that $\varphi(n)/2^{k+1} > 2a+2$ for every $n > c$ when $r \geq 1$, and there is a constant $c$ such that $\varphi(n)/2^k > 2a+2$ for every $n > c$ when $r = 0$.

Now, we prove that for every integer $a \geq 1$, there is always a constant $c$, such that when $n > c$, there exists $x \in Z_n^*$ and $x + 2a \in Z_n^*$ with $x > 1$.

Note that $x_i < n$. If $x_i + 2a > n$ for $2 \leq i \leq t$, then $x_i$ is in the interval $(n - \varphi(n)/2^{k+1} + 2, n)$ when $r \geq 1$ and $x_i$ is in the interval $(n - \varphi(n)/2^k + 2, n)$ when $r = 0$. But there are at most $[\varphi(n)/2^k] - 2$ integers in these two intervals. And it is impossible (since $t - 1 > [\varphi(n)/2^k] - 2$). Therefore there must be some $i$ such that $x \in Z_n^*$ and $x + 2a \in Z_n^*$ with $x > 1$ and Theorem 1 holds.

**Corollary 4:** There exists $x \in Z_n^*$ such that $x > 1$ and $x + 2 \in Z_n^*$ when $n > 6$.

**Remark 1:** By induction, one also can prove that Corollary 4 is true. In [6], we considered the number-theoretic functions $f_1(x) = x$ and $f_2(x) = 2x + 1$, and obtained a similar result as follows:

**Theorem 2:** If $n > 1$ is a positive integer satisfying $n \neq 2, 3, 4, 5, 6$ and $15$, then, there is a positive integer $x > 1$ such that $x \in Z_n^*$ and $2x + 1 \in Z_n^*$, where $Z_n^* = \{x \in N \mid 1 \leq x < n, (x, n) = 1\}$.

**Remark 2:** For the details of proof of Theorem 2, see [**Appendix A**]. However, using the methods in this note, it is difficult to solve the generic case that for any integers $a$ and $b$ with $a > 0$, $(a, b) = 1$ and $2 \mid ab$, there is always a constant $c$, such that when $n > c$, there exists $x \in Z_n^*$ and $ax + b \in Z_n^*$ with $x > 1$. Until recently, we found a simple method which enables us to give slightly stronger results than stated in Corollary 4 and Theorem 2, respectively. Moreover, we obtained an equivalent form of Dickson's conjecture as follows:

Let $s \geq 1$, $f_i(x) = a_i + b_i x$ with $a_i, b_i$ integers for $i = 1, \cdots, s$, if there is a constant $c$ such that for every positive integer $m > c$, there exists a positive integer $x$ such that $f_1(x) > 1$,..., and $f_s(x) > 1$ are all in $Z_m^*$, then there exist infinitely many natural numbers $m$ such that all numbers $f_1(m), \cdots, f_s(m)$ are primes.

**Remark 3:** Theorem 1 and Theorem 2 and the equivalent form of Dickson's conjecture are the first to indicate a common property on some special kinds of number-theoretic functions which likely represent infinitely many primes by rich literatures and a lot of research reports. This common property can be described as follows:

Let $f_1(x_1,\cdots,x_k)$ ,..., $f_s(x_1,\cdots,x_k)$ be $s$ multivariable number-theoretic functions. If $f_1(x_1,\cdots,x_k)$,..., and $f_s(x_1,\cdots,x_k)$ represent simultaneously primes for infinitely many integral points $x=(x_1,\cdots,x_k)$, then there is perhaps a constant $c$ such that for every positive integer $m>c$, there exists an integral point $y=(y_1,\cdots,y_k)$ such that $f_1(y_1,\cdots,y_k)>1$ ,..., and $f_s(y_1,\cdots,y_k)>1$ are all in $Z_m^*$.

**Remark 4:** It has been conjectured that there exists an infinite number of Mersenne primes. In 1964, Gillies [**14**] conjectured that the number of Mersenne primes $\leq x$ is about $\dfrac{2}{\log 2}\log\log x$. In 1983, Wagstaff [**15**] conjectured that the number of Mersenne primes $\leq x$ is about $\dfrac{e^\gamma}{\ln 2}\log\log x$, where $\gamma = 0.5772...$ is Euler's constant. The work of Gillies and Wagstaff showed that the answer is perhaps yes. We find further that the number-theoretic function $2^x-1$ also has the aforementioned property. We have the following Theorem 3. For its proof, see [**Appendix B**].

**Theorem 3:** For every positive integer $n>21$, there is $x>1$ such that $2^x-1$ in $Z_n^*$, where $Z_n^* = \{x \in N \mid 1 \leq x < n, \gcd(x,n)=1\}$.

**Remark 5:** One could prove that for every positive integer $n>10$, there are positive integers of the form $x^2+1$ in $Z_n^*$. But, at present, we have not found an elementary method for proving this result. Of course it is true. Similarly, one could give a sufficient condition of the infinitude of primes of the form $x^2+1$ as follows: Let $n \geq 3$ be a positive integer. If $x$ is the smallest integer such that $x^2+1 \in Z_{n!}^*$, then, $x^2+1$ is a prime number. Namely, this implies that for every positive integer $n>2$, there are primes of the form $x^2+1$ in $Z_{n!}^*$.

**Remark 6:** In 1923, Hardy and Littlewood [**3**] conjectured that the number of primes $x^2+1 \leq C$ is asymptotically equal to $\prod\limits_{p>2}\left(1-\dfrac{1}{p-1}\left(\dfrac{-1}{p}\right)\right)\dfrac{\sqrt{C}}{\log C}$. The problem whether $x^2+1$ represents infinitely many primes goes back to Edmund Landau. In 1912, at International Congress of Mathematicians held in Cambridge, Landau gave four problems which were considered as "unattackable at the present state of science" and named also Landau's problems on primes at present.

The problem above is the first problem of four problems and also is called Landau's first conjecture. In fact, Landau's first conjecture is the special case of Bouniakowsky's conjecture [22]. In 1978, Iwaniec [16] proved that there are infinitely many numbers $x$ such that $x^2+1$ has at most two prime divisors. We do not know that for any $m>3$, whether $1^2+1, 2^2+1,...,m^2+1$ is a W sequence or not [20]. Is this problem equivalent to Landau's first conjecture?

After the great work of Friedlander, John and Iwaniec, Henryk [17] in 1998 and Heath-Brown, D. R. [18] in 2001, maybe, the next goal of this line of research is to prove that Landau's first conjecture is true. It should be interesting to see.

**Remark 7:** One also could obtain an equivalent form of Schinzel-Sierpinski's Conjecture as follows:

Let $f_1(x),\cdots,f_s(x)$ be irreducible univariable polynomials $f_1(x),\cdots,f_s(x)$ with integer coefficients. If there is a constant $c$ such that for every positive integer $m>c$, there exists a positive integer $x$ such that $f_1(x)>1$,..., and $f_s(x)>1$ are all in $Z_m^*$, then there exist infinitely many natural numbers $m$ such that $f_1(m),\cdots,f_s(m)$ are all primes.

By Bateman-Horn's conjecture [21], one can deduce that the necessary condition of Schinzel-Sierpinski's Conjecture [5] holds:

Let $f_1(x),\cdots,f_s(x)$ be irreducible univariable polynomials $f_1(x),\cdots,f_s(x)$ with integer coefficients. If there exist infinitely many natural numbers $m$ such that $f_1(m),\cdots,f_s(m)$ are all primes, then there is a constant $c$ such that for every positive integer $m>c$, there exists a positive integer $x$ such that $f_1(x)>1$,..., and $f_s(x)>1$ are all in $Z_m^*$.

**Remark 8:** However, $2^{2^x}+1$ does not have the aforementioned property. Namely, there is not a constant $c$, such that for every positive integer $m>c$, there is always a positive integer $x$ such that $2^{2^x}+1$ in $Z_m^*$. Does it imply that the number of Fermat primes is finite? Hardy and Wright [7] conjectured that the number of Fermat primes is finite. Nevertheless, they did not give any reasons and explanations. By factoring Fermat number, many people believe that the conjecture in [7] holds. Does our work explain this phenomenon? We will consider further this problem in another paper *On the Infinitude of Some Special Kinds of Primes*.

**Remark 9**: Theorems 1, 2 and 3 are not strong. Using sieve theory [10~13], one can not only obtain the existence of such an $x$ in our theorems, but also some asymptotic formulae of the number of such an $x$. More generally, let $f_1(x_1,\cdots,x_k),\ldots,f_s(x_1,\cdots,x_k)$ be $s$ multivariable number-theoretic functions, one can generalize Euler's totient function as follows:

$$\Phi_{f_1,\ldots,f_s}(n) = \#\{f_1(X) \in Z_n^*,\ldots,f_s(X) \in Z_n^* \mid X = (x_1,\ldots,x_k) \in N^k\}.$$

Obviously, when $s = k = 1$ and $f(x) = x$, $\Phi_{f_1,\ldots,f_s}(n)$ is exactly Euler's totient function $\varphi(n)$. With current technology, it is possible to get some interesting results on the estimation of $\Phi_{f_1,\ldots,f_s}(n)$. But it lies outside the scope of this note and we have not pursued this issue. This should become the subject of future publications. After proving the existence of such $x$ in our theorem, another goal is to give a new sufficient condition that there are an infinite number of twin primes. For details, see Sections 3 and 4.

## 3 A PROOF FOR THE INFINITUDE OF PRIMES

In this section, we begin with Euclid's proof of the infinitude of primes. Euclid's beautiful proof by contradiction goes as follows: Suppose that there are only finitely many primes, say $k$ of them, which denoted by $2 = p_1 < p_2 < \ldots < p_k$. Note that $p_1 p_2 \ldots p_k + 1 > 1$ and hence it must have a prime factor. And this must be $p_j$ for some $j$ with $1 \le j \le k$. But it is impossible since $p_j$ divides both $p_1 p_2 \ldots p_k + 1$ and $p_1 p_2 \ldots p_k$.

Euclid's proof is essentially to construct a number $\beta$ such that $\beta$ is coprime to the product $p_1 p_2 \ldots p_k$. Note that 2 and 3 are prime. So, $|Z_{p_1 \cdots p_k}^*| > 1$ by Euler function formula. On the other hand, as we know, if $a$ is the smallest integer such that $a > 1$ and $(a, p_1 p_2 \ldots p_k) = 1$ then $a$ is prime. Therefore, there are infinitely many primes since $|Z_{p_1 \cdots p_k}^*| > 1$ implies that there is such an integer $a$ in $Z_{p_1 \cdots p_k}^*$. This gives a proof for the infinitude of primes. Although the proof perhaps is not new, it is enlightened us. This proof need not construct a new number $\beta$ such that $\beta$ is coprime to the product $p_1 p_2 \ldots p_k$ but prove directly that there is a number $\beta > 1$ such that $\beta$ is coprime to the product $p_1 p_2 \ldots p_k$. Hence $\beta$ has a new prime factor and it leads to a contradiction. By the existence of such a $\beta$, there must be the least positive integer which is coprime to the product $p_1 p_2 \ldots p_k$. Of course, it is prime.

By Corollary 4, for every positive integer $n > 6$, there is $y > 1$ such that $y$ and $y + 2$ in $Z_n^*$, where $Z_n^* = \{x \in N \mid 1 \leq x < n, \gcd(x,n) = 1\}$. Thus there must be the least positive integer $x > 1$ such that $x$ and $x + 2$ in $Z_n^*$. Note that if $a$ is the smallest integer such that $a > 1$ and $(a, m) = 1$ then $a$ is prime. Naturally, we hope that this property can be preserved. Namely, we hope that such $x$ and $x + 2$ represent simultaneously primes. Unfortunately, this is not true. For example, 7 is the least positive integer which is greater than 1 such that 7 and 9 are in $Z_{10}^*$, but 9 is not primes. How to treat with it? In next section, we will try to give a possible answer.

## 4 A NEW SUFFICIENT CONDITION

Let's look back Euclid's proof again. He considered the product $p_1 p_2 ... p_k$ of primes. Similarly, we may consider $p_k!$. In fact, $p_k! + 1$ and $p_1 p_2 ... p_k$ are coprime, which implies the infinitude of primes again. Directly or more expediently, we consider the factorial $n!$ instead of the finite product $p_1 p_2 ... p_k$ of primes. Clearly, so long as $n > p_k$, then it will lead to a contradiction still. Particularly, let $a \in Z_{n!}^*$ be the smallest integer such that $a > 1$ and $(a, n!) = 1$, then $a$ is prime. This is a key fact. We hope naturally this key fact still is true in the generic case that $f_1(x_1, \cdots, x_k), ..., f_s(x_1, \cdots, x_k)$ represent simultaneously primes for infinitely many integral points. Another reason that we would like to consider the factorial is because the factorial can be viewed as a special case of the $\Gamma$ function which is closely related to the distinguished Riemann Hypothesis. Thus, we naturally hope that for $n > 3$, if $x > 1$ is the smallest integer such that $x$ and $x + 2$ are all in $Z_{n!}^*$, then, $x$ and $x + 2$ are all primes. This is exactly our conjecture in Section 1, which suggests further that it must be infinitely many twin primes. However, we have not been able to work out a complete proof. We tested that the conjecture holds in the many cases for distinct numbers $n$ by using Mathematica. But due to the fact that the value of $n!$ increases rapidly, we can not further test whether it holds for some large integer $n$.

### ACKNOWLEDGEMENTS

Thank Center for Advanced Study in Tsinghua University for providing me with excellent conditions. This work was partially supported by the National Basic Research Program (973) of China (No. 2007CB807902) and the Natural Science Foundation of Shandong Province (No. Y2008G23).

# APPENDIX A

### The Proofs of Theorem 2

**Lemma 5:** Let $p > 5$ be a prime number. For any given integer $a$, there is a positive integer $b$ such that $a(\mod p) \neq b$ and $(b(2b+1)(2a-b)(2(2a-b)+1), p) = 1$.

**Proof:** For any given integer $a$, $b(2b+1)(2a-b)(2(2a-b)+1) \equiv 0(\mod p)$ has at most four solutions in $b$. Note that $\varphi(p) \geq 6$ since $p > 5$. Therefore Lemma 5 is true.

**Lemma 6:** Let $p > 5$ be a prime number and $m$ be a positive integer satisfying $(m, p) = 1$. If there is a positive integer $a$ such that $a, 2a + 1 \in Z_m^*$, then there is a positive integer $x > 1$, such that $x \in Z_{mp}^*$ and $2x + 1 \in Z_{mp}^*$.

**Proof:** Suppose that $a, 2a + 1 \in Z_m^*$, where $Z_m^* = \{x \in N \mid 1 \leq x < m, (x, m) = 1\}$. By Lemma 5, there is a positive integer $b$ such that $a(\bmod\ p) \neq b$ and $(b(2b+1)(2a-b)(2(2a-b)+1), p) = 1$. Suppose that $ml \equiv 1(\bmod\ p)$. Set $y = m((b-a)l(\bmod\ p)) + a$ and $z = m(-(b-a)l(\bmod\ p)) + a$. Clearly, we have

$$(y(2y+1), mp) = 1, (z(2z+1), mp) = 1, 1 < y \leq m(p-1) + a < mp \text{ and } 1 < z < mp.$$

Note that one of $2y + 1$ and $2z + 1$ must be less than $mp$. If $2y + 1 < mp$, then let $x = y$. Otherwise let $x = z$. It shows that Lemma 6 holds.

**Lemma 7:** Let $n > 1$ be a positive integer satisfying $n \neq 3, 6, 15$ and $3 \mid n$, then there is a positive integer $x > 1$ such that $x \in Z_n^*$ and $2x + 1 \in Z_n^*$.

**Proof:** If $p^2 \mid n$, where $p$ is prime, then we may choose $x = -1 + n/p$ such that Lemma 7 holds. Now we suppose that $n > 1$ is a square-free integer. Let $n = 3m$, where $m$ satisfies $(m, 3) = 1$. Since $n \neq 3$, hence $m > 1$ is a square-free integer, too.

(i) If $m$ has only a prime factor $p$, then $p > 5$ since $n \neq 6, 15$. Choose $x = 2$, then $2, 5 \in Z_{3p}^* = Z_n^*$.

(ii) If $m$ has only two prime factors $p$ and $q$. We write $m = pq$ with $p < q$.

 (a) If $(10, pq) = 1$, then choose $x = 2$.

 (b) If $(10, pq) = 2$, then $q \neq 5$. If $q \neq 11$, then set $x = 5$; If $q = 11$, then $n = 66$ and choose $x = 23$.

 (c) If $(10, pq) = 5$, then $p = 5$. If $q \neq 17$, then set $x = 8$. If $q = 17$, then $n = 255$ and choose $x = 11$.

(d) If $(10, pq) = 10$, then $p = 2$, $q = 5$ and $n = 30$ choose $x = 11$.

Thus, one can proved that there always is a positive integer $x \in Z_{3m}^* = Z_n^*$ such that $2x + 1 \in Z_n^*$ when $m$ has only two prime factors. Now we suppose that there is a positive integer $x \in Z_{3m}^*$ such that $2x + 1 \in Z_{3m}^*$ when $m$ has $k \geq 2$ prime factors. We will prove that there is a positive integer $x \in Z_{3m}^*$ such that $2x + 1 \in Z_{3m}^*$ when $m$ has $k + 1 \geq 3$ prime factors. Note that $(m, 3) = 1$. So $m$ has a prime factor $p$ satisfying $p > 5$. Write $m = pr$. We have $y \in Z_{3r}^*$ such that $2y + 1 \in Z_{3r}^*$. By Lemma 6, we know that there is a positive integer $x \in Z_{3rp}^*$ such that $2x + 1 \in Z_{3rp}^*$. Note that $x \neq 1$ (otherwise $2x + 1 = 3$). By induction, when $n = 3m$ is a square-free integer, there is a positive integer $x > 1$ such that $x \in Z_n^*$ and $2x + 1 \in Z_n^*$. Therefore Lemma 7 is true.

**Proof of Theorem 2:** When $n \neq 2, 3, 4, 5, 6$ and $15$, clearly, there is not a positive integer $x > 1$ such that $x \in Z_n^*$ and $2x + 1 \in Z_n^*$. If $3 | n$, clearly, then the theorem is true by Lemma 7. If 3 does not divide $n$ and so does 7, then Theorem is true by setting $x = 3$. Now suppose that 3 does not divide $n$ and set $7 | n$. We write $n = 7^t m$, where $m$ satisfying $(21, m) = 1$ is a positive integer. If $t \geq 2$, then set $x = 7m - 1$. Now we consider the case $t = 1$. If $m = 1$, set $x = 2$. If $m = 2$, set $x = 5$. If $m = 4$, set $x = 5$. If $m = 5$, set $x = 11$. If $m = 6$, set $x = 5$. If $m > 7$, then 3 and $7 \in Z_m^*$, by Lemma 6, we must have $x > 1$ such that $x \in Z_n^*$ and $2x + 1 \in Z_n^*$. So, according to the discussion above, Theorem is true.

Clearly, Theorem 2 also implies that there exists an infinite sequence of positive integers $x_1, ..., x_m, ...$ such that $x_1(2x_1 + 1), ..., x_m(2x_m + 1), ...$ are pairwise relatively prime. Based on this theorem and a similar analysis in this note, one could also give a new sufficient condition that there are an infinite number of Sophie-Germain primes as follows.

**Conjecture 2:** Let $n > 3$ be a positive integer. If $x > 1$ is the smallest integer such that $x \in Z_{n!}^*$ and $2x + 1 \in Z_{n!}^*$, then $x$ and $2x + 1$ must represent simultaneously primes.

# APPENDIX B

## The Proofs of Theorem 3

**Proof of Theorem 3:** Firstly, we generalize the prime-counting function $\pi(x)$ which is the number of primes less than or equal to some real number. Note that pairwise distinct primes are pairwise relatively prime. Consider the number-theoretic function $f(x) = x$ between $N$ and $N$. For any given positive integer $x > 1$, consider a special sub-set $H$ of $\{1, 2, ..., x\}$ as following: $\forall a \in H$, we have $a > 1$, and $\forall a \neq b \in H$, we also have $(a, b) = 1$. Namely, the elements of $H$ are pairwise relatively prime. Denote the set of all such sub-sets of $\{1, 2, \cdots, x\}$ by $\mathrm{M}$. Thus, $\mathrm{M} = \{H \subseteq \{1, 2, ..., x\} \mid \forall a \neq b \in H, (a, b) = 1\}$. Clearly, $\pi(x) = \max_{H \subseteq \mathrm{M}}\{|H|\}$. Namely, $\pi(x)$ can be viewed as the largest among the cardinality of all sub-sets (in which each element exceeds 1 and pairwise distinct elements are pairwise relatively prime) of $\{1, 2, \cdots, x\}$.

Now, let $f(x)$ be a generic number-theoretic function which is monotone increasing. Let $H$ be any sub-set of the image of $f$. Consider the set

$$\mathrm{M} = \{H \subseteq \{1, 2, ..., x\} \mid \forall f(a) \in H, f(a) > 1, \forall f(a) \neq f(b) \in H, (f(a), f(b)) = 1\}.$$

Let $\Pi(x) = \max_{H \subseteq \mathrm{M}}\{|H|\}$. Then, $\Pi(x)$ can be viewed as the generalization of $\pi(x)$. Denote the number of distinct prime factors of $x$ by $\omega(x)$. If we have $\Pi(m) > \omega(m)$, then there is a positive integer $x$ such that $f(x)$ is in $Z_m^*$, and $\Phi_{f_1,...,f_s}(n) > 1$. More generally, let $f_1(x_1, \cdots, x_k), ..., f_s(x_1, \cdots, x_k)$ be $s$ multivariable number-theoretic functions, consider the set

$$\mathrm{M} = \{H \subseteq \{1, 2, ..., x\} \mid 1 \notin H, \forall f_1(a) \neq f_1(b) \in H, ..., f_k(a) \neq f_k(b) \in H, (\prod_{i=1}^{s} f_i(a), \prod_{i=1}^{s} f_i(b)) = 1\},$$

where integral points $a, b$ should be viewed as vectors. Let $\Pi_{f_1,...,f_k}(x) = \max_{H \subseteq \mathrm{M}}\{|H|\}$. Then, $\Pi_{f_1,...,f_k}(x)$ can be viewed as the generalization of $\pi(x)$. Thus, if $f_1(x_1, \cdots, x_k), ...,$ and $f_s(x_1, \cdots, x_k)$ represent simultaneously primes for infinitely many integral points $x = (x_1, \cdots, x_k)$, then, we must have $\Pi_{f_1,...,f_k}(x) \to \infty$ as $x \to \infty$.

Now, let us consider number-theoretic functions $f(x) = 2^x - 1$. Let $2^{r+1} > n \geq 2^r$. Note that $(2^i-1, 2^j-1) = 2^{(i,j)} - 1$. Therefore, it is easy to prove $\Pi(n) = \pi(r)$. But $\pi(r) \geq \dfrac{r}{\log r}$ when $r \geq 17$ [8]. Note also that $\omega(n) \leq \dfrac{\log n}{\log \log n - 1.1714}$ for $n \geq 26$ by Robin's results [9]. Moreover, when $r \geq 109$, we have $\dfrac{1}{4}\log r - 1.1714 > 0$. So, $\log r - 1.1714 > \dfrac{3}{4}\log r$, $\dfrac{r}{\log r} > \dfrac{\frac{3}{4}r}{\log r - 1.1714}$. We also have $\dfrac{3}{4}r > (1+r)\log 2$. Hence, $\dfrac{r}{\log r} > \dfrac{(1+r)\log 2}{\log r - 1.1714} > \dfrac{(1+r)\log 2}{\log r + \log \log 2 - 1.1714}$. Note that $2^{r+1} > n \geq 2^r$. Therefore, we have $(1+r)\log 2 > \log n$ and $\log \log n \geq \log r + \log \log 2$. So, when $r \geq 109$, $\dfrac{r}{\log r} > \dfrac{\log n}{\log \log n - 1.1714}$. Namely, $\Pi(r) = \pi(r) > \omega(n)$ and the theorem is true. When $r < 109$, one could check directly the theorem is true. This completes the proof of the theorem.

**Remark 10:** Using $\sum_{i=1}^{k} \log p_i > p_k - \dfrac{p_k}{2\log p_k}$ in [8], recently, we proved that for every positive integer $k \geq 3$, $2^k \prod_{i=1}^{k} p_i > 2^{p_{k+1}}$ which also implies Theorem 3 holds.

By Theorem 3, for every positive integer $n > 21$, there is $x > 1$ such that $2^x - 1$ in $Z_n^*$. Thus there must be the least positive integer $x > 1$ such that $2^x - 1$ in $Z_n^*$. Of course, $x$ must be prime. But, perhaps, $2^x - 1$ is not prime. For example, $2^{11} - 1 \in Z^*_{(2^2-1)(2^3-1)(2^5-1)(2^7-1)} = Z^*_{82677}$ is the least positive integer of the form $2^x - 1$, but $2^{11} - 1 = 23 \times 89$ is not prime number. Based on the same idea in Sections 3 and 4 of this note, one also could give the following conjecture which suggests that Mersenne primes must be infinitely many.

**Conjecture 3:** Let $n > 3$ be a positive integer. If $x > 1$ is the smallest integer satisfying $2^x - 1 \in Z^*_{n!}$, then, $2^x - 1$ is prime.

# APPENDIX C

## Another proof of Theorem 1

David R Grant and the author further give the following simple proof of theorem 1 [Private Communication].

**Lemma 1:** If an integer $n > 120$, then $n$ has either a prime divisor $p$ such that $p \geq 7$ or a divisor $d$ such that $d \geq b$, where $b \in \{9, 16, 25\}$.

**Proof:** Easy.

**Lemma 2:** Let $p$ be a prime number, and let $a$, $r$ be positive integers. Then, there is a positive integer $y > 1$ such that $\gcd(y(y+2a)(y-6a)(y-4a), p^r) = 1$ for $p \geq 7, r \geq 1$ (resp. $p = 5, r \geq 2$; $p = 3, r \geq 2$; $p = 2, r \geq 4$).

**Proof:** Consider the residue class of $\bmod\ p^r$. $y, y+2a, y-6a, y-4a$ and $y \neq 1 \bmod p^r$ are at most in five residue classes $\bmod\ p^r$ respectively. But, we have that $\varphi(p^r) > 5$ if $p \geq 7, r \geq 1$ (resp. $p = 5, r \geq 2$; $p = 3, r \geq 2$; $p = 2, r \geq 4$). It shows that Lemma 2 holds.

**The Proof of Theorem 1:** Let $n$ be a positive integer $> c$, where $c = \max\{120, 8a+1\}$. By Chinese Remainder Theorem and Lemmas 1 and 2, it is easy to prove that there is a residue class $y \bmod n$ such that $\gcd(y(y+2a)(y-6a)(y-4a), n) = 1$, and $y$ is not $1 \bmod n$. Next, we will prove that there is a positive integer $x$ such that $x > 1$, $x \in Z_n^*$ and $x + 2a \in Z_n^*$. If $y + 2a < n$, then we choose $x = y$. If $y + 2a > n$, then $y + 2a > c \geq 8a + 1$. So, we have that $y > 6a + 1$. Also notice that $1 \neq y < n$. Therefore, we can choose $x = y - 6a$. Thus, we also have that $1 < x < x + 2a = y - 4a < y < n$. It shows that Theorem 1 holds. This completes the proof.

.